\documentclass[a4paper,11pt]{article}
\usepackage{graphicx}
\usepackage{a4wide}
\usepackage[english]{babel}
\usepackage{amsfonts}
\usepackage{amsmath}
\usepackage{amssymb}
\usepackage{dsfont}
\usepackage{amsthm}
\newtheorem{lemma}{Lemma}
\newtheorem{theorem}{Theorem}

\usepackage{float}
\usepackage{url}

\newcommand {\ve} {\varepsilon}

\makeatletter\def\blfootnote{\xdef\@thefnmark{}\@footnotetext}\makeatother
\allowdisplaybreaks
\title{\bf Normal numbers and normality measure}
\author{Christoph Aistleitner\footnote{Department of Applied Mathematics, School of Mathematics and Statistics, University of New South Wales, Sydney NSW 2052, Australia. \mbox{e-mail}: \texttt{aistleitner@math.tugraz.at}. Research supported by a Schr\"odinger scholarship of the Austrian Research
Foundation (FWF).}}
\begin{document}

\date{}
\maketitle

\blfootnote{{\bf MSC 2010:} Primary Classification: 68R15, Secondary Classification: 11K45, 11K16}
\blfootnote{{\bf keywords:} pseudorandom sequence, normality measure, discrepancy, normal numbers}

\begin{abstract}
The normality measure $\mathcal{N}$ has been introduced by Mauduit and S{\'a}rk{\"o}zy in order to describe the pseudorandomness properties of finite binary sequences. Alon, Kohayakawa, Mauduit, Moreira and R{\"o}dl proved that the minimal possible value of the normality measure of an $N$-element binary sequence satisfies
$$
\left( \frac{1}{2} + o(1) \right) \log_2 N \leq \min_{E_N \in \{0,1\}^N} \mathcal{N}(E_N) \leq 3 N^{1/3} (\log N)^{2/3}
$$
for sufficiently large $N$. In the present paper we improve the upper bound to $c (\log N)^2$ for some constant $c$, by this means solving the problem of the asymptotic order of the minimal value of the normality measure up to a logarithmic factor, and disproving a conjecture of Alon \emph{et al.}. The proof is based on relating the normality measure of binary sequences to the discrepancy of normal numbers in base 2.
\end{abstract}

\section{Introduction and statement of results}

Let a finite binary sequence $E_N = (e_1, \dots, e_N) \in \{0,1\}^N$ be given. For $k \geq 1, ~M \geq 1$ and $X \in \{0,1\}^k$, we set 
$$
T(E_N,M,X) = \# \left\{ n:~0 \leq n < M, ~\textrm{and}~(e_{n+1},\dots,e_{n+k}) = X \right\},
$$
which means that $T(E_N,M,X)$ counts the number of occurrences of the pattern $X$ among the first $M+k$ elements of $E_N$. The normality measure $\mathcal{N}(E_N)$ is defined as
\begin{equation} \label{defno}
\mathcal{N} (E_N) = \max_{1 \leq k \leq \log_2 N}~\max_{X \in \{0,1\}^k} ~\max_{1 \leq M \leq N+1-k} \left| T(E_N,M,X) - \frac{M}{2^k}\right|.
\end{equation}
The normality measure has been introduced in 1997 by Mauduit and S{\'a}rk{\"o}zy~\cite{ms1}, together with several other measures of pseudorandomness for finite binary sequences\footnote{Strictly speaking, Mauduit and S{\'a}rk{\"o}zy defined their pseudorandomness measures for sequences over the alphabet $\{-1,1\}$ (instead of $\{0,1\}$, as in the present paper). However, in the case of the normality measure the numerical values of the digits $e_n$ are of no significance whatsoever, since they are used as mere symbols. In the present paper, it is more convenient for our purpose to study sequences defined over the alphabet $\{0,1\}$ (since they can be related to the binary representation of real numbers), and the definitions have been modified accordingly.}. In two papers, Alon, Kohayakawa, Mauduit, Moreira and R{\"o}dl~\cite{akmmr_min,akmmr_typ} studied the \emph{minimal} and the \emph{typical} values of the normality measure (and other measures of pseudorandomness). Concerning the typical value of $\mathcal{N}$, they proved that for any $\ve>0$ there exist $\delta_1, \delta_2 > 0$ such that for $E_N$ uniformly distributed in $\{0,1\}^N$
$$
\delta_1 \sqrt{N} \leq \mathcal{N} (E_N) \leq \delta_2 \sqrt{N}
$$
holds with probability at least $1 - \ve$ for sufficiently large $N$, and conjectured that a limit distribution of
$$
\frac{\mathcal{N} (E_N)}{\sqrt{N}}
$$
exists; the latter has been recently confirmed~\cite{ain}. Concerning the minimal value of $\mathcal{N}$, Alon \emph{et al.} proved that
\begin{equation} \label{eq1}
\left( \frac{1}{2} + o(1) \right) \log_2 N \leq \min_{E_N \in \{0,1\}^N} \mathcal{N}(E_N) \leq 3 N^{1/3} (\log N)^{2/3}
\end{equation}
for sufficiently large $N$. The lower bound in~\eqref{eq1} is based on a relatively simple combinatorial argument. The proof of the upper bound in~\eqref{eq1} is rather elaborate; however, it is entirely constructive, using an explicit algebraic construction based on finite fields. Concerning an possible improvement of~\eqref{eq1}, Alon \emph{et al.} write in~\cite{akmmr_min}
\begin{quote}
``We suspect that the logarithmic lower bound in [equation~\eqref{eq1}] is far from the truth.''
\end{quote}
and formulate the open problem
\begin{quote}
``Is there an absolute constant $\alpha > 0$ for which we have 
$$
\min_{E_N} \mathcal{N}(E_N) > N^\alpha
$$
for all large enough N?''
\end{quote}
In~\cite{akmmr_typ} they write
\begin{quote}
``The authors believe that the answer to [the open problem above] is positive.''
\end{quote}

The purpose of the present paper is to close the gap between the lower and upper bound in~\eqref{eq1}, and settle the problem asking for the asymptotic order of the minimal normality measure of binary sequences, up to a logarithmic factor. More precisely, we will prove that 
\begin{equation} \label{eq2}
\min_{E_N \in \{0,1\}^N} \mathcal{N}(E_N) = \mathcal{O} \left( (\log N)^2 \right),
\end{equation}
by this means giving a negative answer of the problem of Alon \emph{et al.} and disproving their conjecture.

\begin{theorem} \label{th1}
There exists a constant $c$ such that
$$
\min_{E_N \in \{0,1\}^N} \mathcal{N}(E_N) \leq c (\log N)^2
$$
for sufficiently large $N$.
\end{theorem}

The key ingredient of the proof of Theorem~\ref{th1} is to relate the problem asking for binary sequences having small normality measure to the problem asking for normal numbers having small discrepancy. We will describe definitions, basic properties and important results concerning normal numbers in Section~\ref{sec2} below; the proof of Theorem~\ref{th1} will be given subsequently in Section~\ref{sec3}. The proof of Theorem~\ref{th1} is constructive, providing a more or less explicit example of a sequence satisfying the upper bound in the theorem.

\section{Normal numbers} \label{sec2}

Normal numbers have been introduced by Borel~\cite{borel} in 1909. Let $z \in [0,1)$ be a real number, and denote its binary expansion by
$$
z = 0. z_1 z_2 z_3 \dots.
$$
Then $z$ is called a \emph{normal number} (in base 2, which is the only base that we are interested in in the present paper) if for any $k \geq	 1$ and any block of digits $X \in \{0,1\}^k$ the relative asymptotic frequency of the number of appearances of $X$ in the binary expansion of $z$ is $2^{-k}$. Using the terminology from the previous section and writing $Z_N = (z_1, \dots, z_N)$ for the sequence of the first $N$ digits of $z$, this can be expressed as
$$
\lim_{N \to \infty} \frac{T(Z_N,N+1-k,X)}{N} = 2^{-k},
$$
where $k$ is the length of $X$. Borel proved that almost all numbers (in the sense of Lebesgue measure) are normal\footnote{This is the first ever appearance of what we call today the \emph{strong law of large numbers}, for the special case of the i.i.d. system of the Rademacher functions on the unit interval.}. There exist many constructions of normal numbers, the first of them being obtained by concatenating the digital representations of the positive integers (Champernowne~\cite{champ}, 1933), primes (Copeland and Erd\H os~\cite{cope}, 1946) and values of polynomials (Davenport and Erd\H os~\cite{dave}, 1952). Deciding whether a given real number is normal or not is a very difficult problem, and it is unknown whether constants such as $\sqrt{2}, \textup{e}$ and $\pi$ are normal or not (see~\cite{bc}).\\

In an informal way, normal numbers (or the corresponding infinite sequences of digits) are often considered as numbers showing ``random'' behavior (which is justified by the aforementioned theorem of Borel). In fact, different variants of the normality property have been considered as a test for pseudorandomness of (infinite) sequences of digits, for example in the monograph of Knuth~\cite{knu} on \emph{The Art of Computer Programming}, and the normality measure of Mauduit and S{\'a}rk{\"o}zy is a quantitative version of such a pseudorandomness test for the case of a \emph{finite} sequence of digits. For a discussion of the connection between normal numbers, pseudorandomness of (finite) sequences, and pseudorandom number generators, see the book of Knuth and the papers of Mauduit and S{\'a}rk{\"o}zy \emph{On finite pseudorandom binary sequences I-VII}, as well as~\cite{bailey,ni}.\\

To proceed further, we need some notation. A sequence of real numbers $(y_n)_{n \geq 1}$ from the unit interval is called \emph{uniformly distributed modulo one} (u.d. mod 1) if for all intervals $[a,b) \subset [0,1)$ the limit relation
$$
\lim_{N \to \infty} \frac{1}{N} \sum_{n=1}^N \mathds{1}_{[a,b)} (y_n) = b-a
$$
holds. The quality of the uniform distribution of a sequence is measured in terms of the \emph{discrepancy} $D_N$, which for $N \geq 1$ is defined as
$$
D_N(y_1, \dots, y_N) = \sup_{0 \leq a < b \leq 1} \left| \frac{1}{N} \sum_{n=1}^N \mathds{1}_{[a,b)} (y_n) - (b-a) \right|.
$$ 
A sequence is u.d. mod 1 if and only if its discrepancy tends to zero as $N \to \infty$.\\

By an observation of Wall~\cite{wall}, a number $z$ is normal (in base 2) if and only if the sequence
$$
\left( \langle 2^{n-1} z\rangle \right)_{n \geq 1},
$$
where $\langle \cdot \rangle$ denotes the fractional part, is u.d. mod 1. Equivalently, $z$ is normal if and only if
$$
D_N \left( z, \langle 2 z \rangle, \dots, \langle 2^{N-1} z \rangle \right) \to 0 \qquad \textrm{as} \qquad N \to \infty.
$$
Korobov~\cite{koro} posed the problem of finding a function $\psi(N)$ with maximal decay for which there exists a number $z$ such that
$$
D_N \left( z, \langle 2 z \rangle, \dots, \langle 2^{N-1} z \rangle \right) \leq \psi(N), \qquad N \geq 1.
$$
The best results concerning this question is currently due to Levin~\cite{levin}, who proved (constructively, by giving an explicit example) the existence of a $z$ for which
\begin{equation} \label{levin}
D_N \left( z, \langle 2 z \rangle, \dots, \langle 2^{N-1} z \rangle \right) = \mathcal{O} \left(\frac{(\log N)^2}{N^{-1}} \right) \qquad \textrm{as} \qquad N \to \infty.
\end{equation}
This result should be compared with a lower bound of Schmidt~\cite{schmidt}, stating that for \emph{any} sequence $(y_n)_{n \geq 1}$
$$
D_N(y_1, \dots, y_N) \geq c_{\textup{abs}} \frac{\log N}{N}.
$$
Thus Korobov's problem is solved, up to a logarithmic factor. It is also interesting to compare~\eqref{levin} with the ``typical'' discrepancy of a normal number: for almost all $z \in [0,1)$,
$$
\limsup_{N \to \infty} \frac{\sqrt{N} D_N \left( z, \langle 2 z \rangle, \dots, \langle 2^{N-1} z \rangle \right)}{\sqrt{\log \log N}} = \frac{2 \sqrt{21}}{9}
$$
(Fukuyama~\cite{fuku}).\\

For more information on normal numbers we refer to~\cite{harm,qu}, for an introduction to uniform distribution and discrepancy theory to~\cite{dts,kn}.\\

The main tool in the proof of Theorem~\ref{th1} is the following lemma.

\begin{lemma} \label{lemma1}
Let $z \in [0,1)$ be a real number, whose binary expansion is given by
$$
z = 0. z_1 z_2 z_3 \dots,
$$
and assume that there exists a nondecreasing function $\Phi(N)$ such that
\begin{equation}
D_N \left( z, \langle 2 z \rangle, \dots, \langle 2^{N-1} z \rangle \right) \leq \frac{\Phi(N)}{N}, \qquad N \geq 1.
\end{equation}
Then for each $N \geq 1$ the binary sequence $Z_N = (z_1, \dots, z_N)$ satisfies
$$
\mathcal{N}(Z_N) \leq \Phi(N).
$$
\end{lemma}

In view of Levin's result~\eqref{levin}, Theorem~\ref{th1} is a direct consequence of the lemma.

\section{Proof of Theorem~\ref{th1}.} \label{sec3}

By the previous remark, to establish Theorem~\ref{th1} it remains to prove Lemma~\ref{lemma1}. Let a number $N$ be fixed, and assume that 
$$
D_N \left( z, \langle 2 z \rangle, \dots, \langle 2^{M-1} z \rangle \right) \leq \frac{\Phi(M)}{M}
$$
holds for $1 \leq M \leq N$. To prove $\mathcal{N}(Z_N) \leq \Phi(N)$ we have to show that for \emph{any} values of $k,~X$ and $M$ satisfying $1 \leq k \leq \log_2 N$, $X \in \{0,1\}^k$ and $1 \leq M \leq N-k+1$ we have
\begin{equation} \label{meq}
\left| T(Z_N,M,X) - \frac{M}{2^k}\right| \leq \Phi(N).
\end{equation}
Let $k, ~X$ and $M$ satisfying these assumptions be fixed and write $X = (x_1, \dots, x_k)$, where $x_1, \dots, x_k \in \{0,1\}$. By definition,
\begin{eqnarray*}
T(Z_N,M,X) & = & \# \left\{ n:~0 \leq n < M, ~\textrm{and}~(z_{n+1},\dots,z_{n+k}) = (x_1, \dots, x_k) \right\}.
\end{eqnarray*}
To $X$ we can assign an interval $I_X$ by setting
$$
I_X = \left[ \sum_{j=1}^k x_j 2^{-j}, \left(\sum_{j=1}^k x_j 2^{-j}\right) + 2^{-k} \right).
$$
Then $I_X$ is a half-open interval of length $2^{-k}$. The following observation is the crucial point of the proof of the lemma. We have
$$
(z_{n+1},\dots,z_{n+k}) = (x_1, \dots, x_k)
$$
if and only if
$$
\langle 2^{n} z \rangle \in I_X.
$$
In fact, we have
$$
\langle 2^{n} z \rangle = 0.z_{n+1} z_{n+2} \dots,
$$
and for any number $y \in [0,1)$ the relation $y \in I_X$ holds if and only if the first $k$ digits of $y$ coincide with $(x_1, \dots, x_k)$.\\

Consequently, we have
\begin{eqnarray} \label{u2}
T(Z_N,M,X) & = & \sum_{n=0}^{M-1} \mathds{1}_{I_X} \left( \langle 2^n z \rangle \right) = \sum_{n=1}^{M} \mathds{1}_{I_X} \left( \langle 2^{n-1} z \rangle \right).
\end{eqnarray}
Now by the assumption on the discrepancy of $z$ we have
\begin{equation} \label{u1}
\left|\frac{1}{M} \sum_{n=1}^{M} \mathds{1}_{I_X} \left( \langle 2^{n-1} z \rangle \right) - \frac{1}{2^k} \right| \leq \frac{\Phi(M)}{M},
\end{equation}
and consequently, multiplying equation~\eqref{u1} by $M$ and using~\eqref{u2}, we obtain
$$
\left| T(Z_N,M,X) - \frac{M}{2^k} \right| \leq \Phi(M).
$$
Since by assumption the function $\Phi(M)$ is nondecreasing, this establishes~\eqref{meq}, which proves Lemma~\ref{lemma1}.

\end{document}